\documentclass[a4paper,twoside,fleqn,12pt]{article}

\usepackage{epsfig}
\usepackage[shortlabels]{enumitem}
\usepackage{amsmath}
\usepackage{amssymb}
\usepackage{amsthm}
\usepackage{datetime}

\newcommand{\commentaar}[1]{{}}

\newcommand{\diag}{\text{diag}}
\newcommand{\field}[1]{\mathbb{#1}}
\newcommand{\inprod}[2]{\langle #1 , #2 \rangle}
\newcommand{\partieel}[2]{\frac{\partial #2}{\partial #1}}
\newcommand{\partiel}[2]{\frac{\partial #1}{\partial #2}}

\newcommand{\norm}[1]{\lVert #1 \rVert}

\newcommand{\Lgr}{\mathcal{L}}

\newcommand{\R}{\field{R}}
\newcommand{\SE}{\text{SE}}
\newcommand{\SO}{\text{SO}}

\newcommand{\cF}{{\mathcal F}}
\newcommand{\cM}{{\mathcal M}}

\newcommand{\e}{\textrm{e}}

\DeclareMathOperator{\grad}{grad}

\newcommand{\smallo}{\textrm{o}}

\newtheorem{theorem}{Theorem}[section]
\newtheorem{lemma}[theorem]{Lemma}
\newtheorem{proposition}[theorem]{Proposition}
\newtheorem{corollary}[theorem]{Corollary}
\newtheorem{definition}[theorem]{Definition}
\newtheorem{Remark}[theorem]{Remark}

\newenvironment{remark}{\begin{Remark} \begin{rm}}{\end{rm} \end{Remark}}

\title{Critical points at infinity in charged $N$-body systems}
\author{I. Hoveijn, H. Waalkens, M. Zaman\\
\small{Bernoulli Institute for Mathematics, Computer Science and Artificial Intelligence}\\
\small{University of Groningen, The Netherlands}}

\begin{document}

\maketitle

\section*{Abstract}
We define the notion of critical points at infinity for the charged $N$-body problem, following the approach of Albouy \cite{alb1993}. We give a characterisation of such points and show how they can be found in the charged $3$-body problem. The symmetry group of the $N$-body problem and accompanying integrals play a key role. In fact critical points at infinity are indispensible in understanding the bifurcations of the integral map. Together with the critical points at infinity in the charged $3$-body problem, we present the bifurcation values.

\textbf{Keywords:} Charged $N$-body problem, integral map, critical sequence, cluster decomposition, critical point at infinity, critical point, relative equilibrium

%\tableofcontents

%----------------------------------------------------------------------
\section{Introduction}\label{sec:intro}

%\commentaar{waar gaat het over?}
This is the second article in a series of three on the charged three body problem. We continue our study of the integral manifolds of the charged three body problem, started in \cite{hwz2019}. In the latter we in particular considered the critical points of the integral map. Here we extend this to critical points at infinity of the charged $N$-body problem. In the last section we return to the charged three body problem. In the third article \cite{hwz2019c}, we study the Hill regions of such systems.

We consider the charged $N$ body problem as a Hamiltonian dynamical system, in particular as a mechanical dynamical system (terminology of Arnol'd, see \cite{arnold1980}). That is, the Hamiltonian is the sum of 'kinetic energy' and 'potential energy'. The $N$ bodies reside in $\R^3$, therefore the phase space $\cM$ is the $6N$-dimensional cotangent bundle $\cM = T^* \R^{3N} \simeq \R^{3N} \times \R^{3N}$ of the $3N$-dimensional \emph{configuration space} $\R^{3N}$. The interaction of the bodies is governed by a potential function $V$. Furthermore, the function $V$ is assumed to be $\SE(3)$-invariant. In fact we assume pairwise interaction and in particular that the potential energy of two bodies is inverse proportional to the distance of the two bodies. For further details about the phase space and the Hamiltonian function, see section \ref{sec:set}.

In this setting, the $N$-body problem admits ten integrals: the Hamiltonian function, together with three components of momentum, angular momentum and the centre of mass. They define the so-called \emph{integral map} $\cF: \R^{6N} \to \R^{10}$. In the end we wish to characterise the \emph{integral manifolds} defined as the fibres $\cF^{-1}(\phi)$ of the map $\cF$ for all values $\phi$. Here we first consider critical values of $\cF$ from critical points of $\cF$, for at these values we expect the integral manifolds to change topology. Due to the nature of the potential, the integral manifolds are not necessarily compact for finite values of the integral map. Thus we may expect critical points at infinity.

%\commentaar{waarom van belang?}
Critical points at infinity cannot be found by standard means, but their presence can be conjectured from topological changes in the fibres of $\cF$ at seemingly regular values of $\cF$. This has happened in the three body problem with gravitation, see for example \cite{simo1975}. The first definition of critical points at infinity for $N$-body problems that we are aware of is by Albouy, see \cite{alb1993}. Indeed, the critical values corresponding to critical points at infinity in the three body problem with gravitation are those referred to in Simo, see \cite{simo1975}. Needless to say that critical points and critical points at infinity are of equal importance in understanding the topology of the integral manifolds.

%\commentaar{wat is essentieel?}
Three properties of the $N$-body problem are essential for this study. First the presence of the symmetry group $\SE(3)$. This provides us with the integral map and thus integral manifolds. Another implication, together with the assumption of pairwise interaction, is that the potential $V$ depends on inter body distances only. The second property is that the potential and its derivatives go to zero for inter body distance going to infinity. This implies that there is no interaction between bodies going infinitely far apart. For instance in a three body problem we have the possibility of one body at rest and infinitely far away two bodies forming a two body problem of bodies rotating around their centre of mass. This will be an example of critical points at infinity (in fact a critical orbit). The third property is the homogeneity of the potential function. This property implies a non-symplectic scaling symmetry of the Hamiltonian system. The Hamiltonian is not invariant under this scaling, but with proper time scaling the equations of motion are invariant. From this property we get families of orbits and so families of critical points. These three properties will be treated in more detail in section \ref{sec:set}.

%\commentaar{wat is het resultaat?}
In order to give a precise account of the results on critical points at infinity we would need the notions of section \ref{sec:scc}. Therefore we state these results here in an informal way, for a precise statement see section \ref{sec:cpimain}. Suppose that the $N$-body problem is split into $k_i$-body subsystems such that $\sum_i k_i = N$. Furthermore suppose these $k_i$-body subsystems are mutually infinitely far apart in configuration space. Loosely speaking the critical points at infinity of the $N$-body problem consist of the critical points of these subsystems. In general the number of critical points of a $k$-body problem is unknown, see for example \cite{acs2012}. But in case of three body problems, the subsystems consist of one or two bodies only. It turns out that in a three body problem with Coulomb interaction, with positive and negative charges, we find two critical points at infinity. Thus we find at most two critical values of the integral map related to critical points at infinity. In the gravitational case we find three critical points at infinity with at most three critical values, see \cite{alb1993}.

%\commentaar{hoe zit dit artikel in elkaar, hoe gaan we te werk}
We proceed as follows. In section \ref{sec:set} we specify the charged $N$-body problem as a mechanical dynamical system. Here we also detail several properties from general to more specific. All such systems will be symplectic and moreover $\SE(3)$-invariant. The latter allows for a number of integrals. As soon as the potential is defined we see that for our choice the system also admits a non-symplectic scaling symmetry.

We are very much indebted to the work of Albouy on the $N$-body problem with gravitation, see \cite{alb1993}. In fact we do not claim originality for the constructions we describe in sections \ref{sec:scc} and \ref{sec:promain}. We adapted the constructions from \cite{alb1993} with slight extensions, generalisations and rearrangements to suit our needs.

In section \ref{sec:scc} a first step towards a definition of a critical point at infinity is by requiring a Lagrange multiplier condition in a limit. This however is not enough to yield a critical point at infinity under all circumstances. So part of this section is devoted to filtering out the desired points by additional conditions. A particular choice of coordinates is very helpful in this respect especially in combination with the notion of clusters. The latter are $l$-body subsystems of the $N$-body problem, whose centres of mass are infinitly far apart, but where the distances of the bodies to the center of mass remain finite within each cluster.

The main results are stated in section \ref{sec:cpimain}. Section \ref{sec:promain} is devoted to the proofs of the results. The last section discusses the meaning of the results in section \ref{sec:cpimain} for the charged $3$-body problem. For more discussion we refer to the third article \cite {hwz2019c} in this series.

%--------------------------------------------------------------------------------
\section{Setting of the problem}\label{sec:set}
We study the $N$-body problem as a symplectic dynamical system. Since each body resides in $\R^3$, the \emph{configuration space} is $\R^{3N}$ and the \emph{phase space} $\cM$ is the cotangent bundle of the configuration space: $\cM = T^* \R^{3N} \simeq \R^{3N} \times \R^{3N}$. Let $\Pi: \cM \to \R^{3N}$ be the projection from phase space to configuration space. On $\cM$ we take the standard \emph{symplectic form} $\omega = \sum_i dp_i \wedge dq_i$ expressed in coordinates $q=(q_1,\ldots,q_N)$ and $p=(p_1,\ldots,p_n)$ with $q_i,p_i\in\R^3$. The Hamiltonian of this system has a special form: $H(q, p) = T(p) + V(q)$, where the physical interpretation of $T$ and $V$ is kinetic energy and potential energy, respectively. Kinetic energy is given by
\begin{equation}\label{eqn:kinegy}
T(p) = \frac{1}{2} \sum_i \frac{p_i^2}{m_i}
\end{equation}
where $p_i$ is momentum of body $i$ and $m_i$ denotes its mass. Potential energy will be given shortly, it is determined by the interaction of the bodies. Although we will not explicitly study dynamics in this article, knowing the equations of motion is still useful. They take the familiar form ($\omega$ in standard form)
\begin{equation}\label{eqn:hameq}
\dot{q}_i = \partiel{H}{p_i} \;\text{ and }\; \dot{p}_i = -\partiel{H}{q_i}.
\end{equation}

From physics we also know that the interaction between bodies is in many cases (at least for gravity and Coulomb interaction) independent of position and orientation of the system in space. Such a system admits $\SE(3)$-symmetry. Here we just restrict to systems with this symmetry. This group acts by translations and rotations on phase space as
\begin{equation}\label{eqn:sedrie}
\begin{aligned}
\tau_a(q_i, p_i) &= (q_i + a, p_i)\\
\rho_g(q_i, p_i) &= (gq_i, gp_i)
%
%\big((q_1, p_1), \ldots, (q_N, p_N)\big) &\mapsto \big(\tau_a(q_1, p_1), \ldots \tau_a(q_N, p_N)\big)\\
%\big((q_1, p_1), \ldots, (q_N, p_N)\big) &\mapsto \big(\rho_g(q_1, p_1), \ldots \rho_g(q_N, p_N)\big)
\end{aligned}
\end{equation}
for $i \in \{1,\ldots,N\}$, $a \in \R^3$ and $g \in \SO(3)$.

Two consequences of the fact that the system is $\SE(3)$-symmetric are relevant for us. The first is the existence of integrals, see for example \cite{arnold1980}, \cite{am1987} or \cite{or2004} for more details on the relation between symmetry groups and integrals. The second is the existence of (generators of) group invariants, see \cite{weyl1997} for the particular case of $\SE(3)$. As a third consequence we could say that the functions $T$ and $V$ are $\SE(3)$-invariant, but this fact is implicit in the definition of the $N$-body problem admitting $\SE(3)$ as a symmetry group.

The integrals related to the translation group are the three components of total momentum $P(q, p) = \sum_i p_i$. The component functions of $P$ have no critical points on phase space, so without loss of generality we may fix $P$ at any value. From now on we assume that $P = 0$, this leaves us with three more integrals. Indeed, from the equations of motion, together with $T(p) = \frac{1}{2} \sum_i \frac{p_i^2}{m_i}$, we immediately see that the components of centre of mass $Q(q, p) = \frac{1}{M} \sum m_i q_i$ where $M = \sum m_i$, are also integrals on the zero level set of $P$. Just like $P$, the component functions of $Q$ have no critical points on phase space and we may fix $Q$ at any value. Integrals related to rotation symmetry are the components of total angular momentum $L(q,p) = \sum_i q_i \times p_i$. These functions do have critical points on phase space and the critical value equals zero. However, we will not consider this value of angular momentum in the sequel. From now on we do assume that the value of $P$ is zero, so that we have the \emph{integral map}
\begin{equation}\label{eqn:intmap}
\cF: \R^{3N} \times \R^{3N} \to \R^{10} : (q,p) \mapsto \big( H(q,p), L(q,p), P(q,p), Q(q,p) \big),
\end{equation}
which we will consider on the zero level set of $P$ and $Q$.

\commentaar{niet zo over die potentiaal zeuren, gewoon 1/r, anders formuleren en korter}

$\SE(3)$-symmetry, by a theorem of Po\'enaru, see \cite{poe1976}, yields that the potential function $V$ is a function of generators of $\SE(3)$-invariants, see \cite{weyl1997}. In this case these are invariants of the action of $\SE(3)$ on configuration space. Now we consider pairwise interaction of the bodies such that potential energy is additive over pairs. In fact we take a potential that can be expressed using a single function $f : \R \to \R$, namely
\begin{equation}\label{eqn:potegy}
V(q) = \sum_{i < j} \alpha_{ij} f\big(\norm{q_i - q_j}\big)
\end{equation}
where the $\alpha_{ij}$ are coefficients determined by properties of the bodies. We will always assume that $f$ and $f'$ tend to zero at infinity. In case of gravitation, $\alpha_{ij} = -m_i m_j$ and in case of electrostatic interaction, $\alpha_{ij} = c_i c_j$ where $c_i$ is the \emph{electric charge} of body $i$. In both cases $f(x) = \frac{1}{x}$. The latter suggests that we may wish to exclude the \emph{collision set}
\begin{equation*}%\label{eqn:colset}
\Delta = \{(q_1,\ldots,q_N) \in \R^{3N} \;|\; q_i = q_j \text{ for } i \neq j \}
\end{equation*}
from the configuration space. Furthermore we assume that the Hamiltonian system is smooth on the cotangent bundle of the configuration space without the collision set.

In the setting of the $\SE(3)$-symmetric $N$-body problem with pairwise interaction, the potential function $V$ is limited to the form in equation \eqref{eqn:potegy}, but still with the possibility to choose the function $f$. In the main examples (gravitation and electrostatic interaction) the function $f$ is singular at zero. But even if $f$ is non-singular it is still interesting to ask whether there may be collisions of bodies. Therefore the notion of collision set is not just a technical detail.

In case of gravitation or electrostatic interaction, the function $f$ and thus potential $V$, is homogeneous of degree $-1$. Generalising this to $f$ being homogeneous of degree $\mu$ we see that the system admits a non-symplectic \emph{dilation symmetry} acting as
\begin{equation}\label{eqn:dilation}
\delta_s(q_i, p_i) = (\e^{2s} q_i, \e^{\mu s} p_i)
\end{equation}
for $i \in \{1,\ldots,N\}$ and $s \in \R$. Although this group does not act symplectically, orbits of the $N$-body problem are taken to orbits, since with a proper scaling of time the equations of motion \eqref{eqn:hameq} are invariant. Indeed, if $(q, p, t) \mapsto \big(\e^{2s} q, \e^{\mu s} p, \e^{(2-\mu)s} t\big)$, the equations of motion are invariant. The Hamiltonian, however, is not invariant: $H(\delta_s(q_i, p_i)) = \e^{2\mu s} H(q_i, p_i)$.

From the considerations above we see that the charged $N$-body problem admits a large symmetry group. Therefore the system can be simplified by applying reduction. However explicit reduction with respect to rotation or dilation symmetry is not going to present much insight at the moment. Reduction with respect to translation symmetry is rather straightforward. In the next section we will introduce coordinates on the translation reduced space, we will call these \emph{Albouy coordinates}. In fact the \emph{translation reduced phase space} can be idenfified with the zero level set of $P$ and $Q$ giving the reduced space $\R^{6N-6}$ and the \emph{reduced integral map}
\begin{equation}\label{eqn:intmapred}
\cF_r: \R^{6N-6} \to \R^4: (q, p) \mapsto \big(H(q, p), L(q, p) \big)
\end{equation}
which will be used in the next sections. The remaining dilation-rotation symmetry generates equivalence classes of critical points of the reduced integral map. These, in turn, generate equivalent critical values of $\cF_r$. Let us look at the induced action on the image of $\cF_r$. If $(h, l)$ is a value of $\cF_r$ at $(q,p)$, then $\cF_r(\rho_R(q, p)) = (h, Rl)$ and $\cF_r(\delta_s(q, p)) = (\e^{-2s} h, \e^s l)$. Thus we see that $h \norm{l}^2$ is constant on equivalent values of $\cF_r$. In fact it is the only generator of the invariants of the induced action. Therefore the scalar $-h \norm{l}^2$ is sometimes called the \emph{bifurcation parameter} of the reduced integral map, see for example \cite{mmw1998}.

%--------------------------------------------------------------------------------
\section{Sequences, coordinates and clusters}\label{sec:scc}
In this section we use sequences satisfying a Lagrange multiplier equation in the limit, to define critical points at infinity. Furthermore we define Albouy coordinates on a translation reduced space. With help of these coordinates we define so called clusters of subsystems which we encounter when studying critical points at infinity. We closely follow Albouy, see \cite{alb1993}, with a few adaptations to our setting.

%--------------------------------------------------------------------------------
\subsection{Sequences}\label{sec:seq}
A sequence of states is simply a sequence $z_k$ of points in phase space $\cM = \R^{6N}$. In the sequel we will need a notion of convergence. To define this we use the norm $\norm{z}$ of $z \in \cM$. We do allow a sequence to converge to infinity and to avoid descriptions like 'finite limit' and 'infinite limit' we introduce the following language.
\begin{definition}[Convergence]\label{def:limits}
We say that a sequence $z_k$ converges if either of the two holds:
\begin{enumerate}[i), itemsep=0.5ex, topsep=0pt]
\item $z_k$ converges to a point $z$ if $\lim_{k \to \infty} \norm{z_k - z} = 0$
\item $z_k$ converges to infinity if $\lim_{k \to \infty} \norm{z_k} = \infty$.
\end{enumerate}
Both in the usual sense as real numbers. Similarly, if $F$ is a real valued function on phase space, we say that the sequence $F(z_k)$ converges if it converges to a real value or to infinity like in the above.
\end{definition}
The motivation to consider sequences is that we wish to define critical points at infinity. Critical points of the integral map $\cF$ are points such that the rank of the Jacobi matrix of $\cF$ is not maximal. Or, put differently, the gradient of the Hamiltonian and the gradients of the integrals $F_i$ are linearly dependent. Using this we are able to find finite critical points. To define critical points at infinity we first define sequences of states asymptotically having the property of linear dependence of the gradients. It turns out that this condition is not sufficient. That is, we also find critical points, but possibly collisions as well. To define critical points at infinity we have to filter out these possibilities. So, first we define sequences having properties that make them converge to (critical) points. Then we define further properties to be able to distinguish sequences converging to (critical) points as opposed to sequences converging to (critical) points at infinity. Finally we define a criterion to filter out sequences converging to collisions.
Suppose the system at hand has $n+1$ independent integrals so that $\cF = (H, F_1,\ldots,F_n)$. The following definition aims to catch both critical points and critical points at infinity by requiring that linear dependence holds in a limit in the sense defined above.
\begin{definition}[Compatible sequence of multipliers of $H$]\label{def:compmul}
Let $z_k$ be a sequence of states. Then $\lambda_k = (\lambda_{1,k}, \ldots, \lambda_{n,k})$ is called a \emph{compatible sequence of multipliers for $H$} if:
\begin{equation*}
\lim_{k \to \infty} \Lgr(z_k, \lambda_k) = 0,
\end{equation*}
where
\begin{equation*}
\Lgr(z, \lambda) =  \grad H(z) - \sum_{i=1}^n \lambda_i\, \grad F_i(z).
\end{equation*}
\end{definition}
Clearly not every sequence of states will have an accompanying sequence of compatible multipliers, but if it has we possibly approach a critical point (finite or infinite).
\begin{definition}[Critical sequence]\label{def:critseq}
If a sequence of states $z_k$ has an accompanying sequence of compatible multipliers for $H$ we call it a \emph{critical sequence of $H$}.
\end{definition}
In the above we implicitly defined a critical point of the integral map $\cF$ as a critical point of $H$ restricted to the levels of $F_1,\ldots,F_n$, which are assumed to be regular. If $H$ has such a critical point then there are (uncountable) many critical sequences converging to that point. Among these there are still uncountable many for which $F_i(z_k),\ldots,F_n(z_k)$ remain constant. If on the other hand we are given a critical sequence $z_k$ of $H$, then it may converge or not. If it converges to a point we conjecture that it either converges to a collision or to a critical point. These possibilities can be distinguished by the limit value of the potential $V$. Here we assume that the function $f$, see section \ref{sec:set}, has the property $|f(x)| \to \infty$ when $x \to 0$. If $V_k$ converges to $-\infty$ the limit of $z_k$ is a collision otherwise it is a critical point. The latter is in fact a consequence of corollary \ref{cor:nocol} (see section \ref{sec:techs}), which says that a (subsequence of) a citical sequence never approaches a collision of repelling bodies. If the critical sequence converges to infinity we consider the limit as a critical point at infinity if the values of the integrals remain finite. Let these arguments suffice to justify the following definition.
\begin{definition}[Horizontal critical sequence]\label{def:horcritseq}
A critical sequence $z_k$ is called \emph{horizontal} if it has the properties:
\begin{enumerate}[i), itemsep=0.5ex, topsep=0pt]
\item $F_1(z_k),\ldots,F_n(z_k)$ are finite and constant, that is independent of $k$
\item $H(z_k)$ converges to a finite value for $k \to \infty$
\end{enumerate}
\end{definition}
With the previous definitions we may also find the finite critical points. Here we are explicitly looking for critical points at infinity. So we exclude the former and finally define critical points at infinity.
\begin{definition}[Critical point at infinity]\label{def:cpinf}
If the projection on configuration space of a horizontal critical sequence $z_k$ does not have a subsequence converging to a point we say that $z_k$ converges to a \emph{critical point at infinity}.
\end{definition}

\begin{remark}\label{rem:notationa}
In the remaining part of the article we will often consider functions like kinetic energy, potential energy, moment of inertia and angular momentum on sequences. If confusion is not likely we write $F_k$ instead of $F(z_k)$ to indicate a sequence of values of the function $F$ on the sequence $z_k$.
\end{remark}

%--------------------------------------------------------------------------------
\subsection{Coordinates}\label{sec:albco}
Following Albouy \cite{alb1993} we define coordinates on a translation reduced space. Since the Hamiltonian is invariant under the action action of the translation group, by Noether's theorem, total momentum is a conserved function on phase space. Then the centre of mass is also a conserved function on phase space, provided that total momentum equals zero, see section \ref{sec:set}. The expression for centre of mass motivates the following definition.
\begin{definition}[Albouy coordinates]\label{def:dispo}
Let $m_1,\cdots,m_N$ be the masses of the bodies in a $N$-body problem, then we define the linear subspace of $\R^N$, $D_N = \{(\xi_1, \ldots, \xi_N) \;|\; \sum_i m_i \xi_i = 0 \}$. Furthermore let $(\cdot,\cdot)$ be the standard inner product on $\R^N$. Define a new inner product $\inprod{\cdot}{\cdot}$ on $\R^N$ as $\inprod{X}{Y} = (X, MY)$ for all $X, Y \in \R^N$ and $M = \diag(m_1,\ldots,m_N)$. The associated norm is $\norm{X}^2 = \langle X, X \rangle$ and $D_N$ inherits this inner product.
\end{definition}
Using this definition we are going to define reduced coordinates on configuration space. Now we switch to velocities instead of momenta and we define the state of a system of $N$ bodies in $\R^3$ as the six-tuple $S = (X, Y, Z, P, Q, R)$ for $X,\ldots,R \in D_N$. Here $X$ contains the first components of the \emph{positions} of the bodies and similarly $Y$ and $Z$ contain the second and third components. The $P$, $Q$ and $R$ in a similar way contain the components of the \emph{velocities} of the bodies. In the phase space $D_N^6$ we again have a \emph{collision set} and sometimes we may wish to exclude these states for obvious reasons.
\begin{definition}[Collision set]\label{def:colset}
Let $X = (\xi_1,\ldots,\xi_N) \in D_N$, $Y = (\eta_1,\ldots,\eta_N) \in D_N$ and $Z = (\zeta_1,\ldots,\zeta_N) \in D_N$, then the set $\Delta = \{(X, Y, Z, P, Q, R) \in D_N^6 \;|\; \xi_i = \xi_j, \eta_i = \eta_j, \zeta_i = \zeta_j \;\text{for}\; i \neq j\}$ is called the \emph{collision set}. The projection of $\Delta$ on configuration space will also be called collision set.
\end{definition}
In the sequel we call coordinates $S$, \emph{Albouy coordinates} on the translation reduced space. The advantages of these coordinates are:
\begin{enumerate}[i), itemsep=0.5ex, topsep=0pt]
\item reduced coordinates without the need to make choices: all bodies of the $N$-body problem are treated equally,
\item kinetic energy, angular momentum and moment of inertia can be concisely expressed, in particular in combination with the notion of \emph{clusters}, see below.
\end{enumerate}
The space $D_N^6$ can be considered as the translation reduced phase space. On this space frequently used functions like potential energy $V$, kinetic energy $T$, total energy $H$, angular momentum $L$ and moment of inertia $I$ have simple expressions.
\begin{lemma}\label{lem:functions} The functions $K$, $H$, $L$ and $I$ take a simple form in Albouy coordinates:
\begin{enumerate}[i), itemsep=0.5ex, topsep=0pt]
\item \emph{Kinetic energy} is given by $\frac{1}{2}K(S)$ where $K(S) = \norm{P}^2 + \norm{Q}^2 + \norm{R}^2$.
\item The \emph{total energy} of the system is $H(S) = \frac{1}{2} K(S) + V(S)$.
\item The three components of \emph{angular momentum} are\\
$L(S) = \big(\inprod{Y}{R} - \inprod{Z}{Q}, \;\inprod{Z}{P} - \inprod{X}{R}, \;\inprod{X}{Q} - \inprod{Y}{P}\big)$.
\item Let $I(S) = \norm{X}^2 + \norm{Y}^2 + \norm{Z}^2$, then the \emph{moment of inertia} is $I(S)$.
\end{enumerate}
\end{lemma}

We still have the action of the rotation group $\SO(3)$. However, elements of $\SO(3)$ do not act in a nice way on coordinates $S$. This can be improved by using the bijection $\phi : \R^{6N} \to \R^{3\times 2N}$ taking $(X, Y, Z, P, Q, R)$ to
\begin{equation*}
\begin{pmatrix} X & P \\ Y & Q \\ Z & R \end{pmatrix}
\end{equation*}
and restrict $\phi$ to $D_N^6$. Then $\phi$ is a bijection $D_N^6 \to \phi(D_N^6)$. Since the action of $g \in \SO(3)$ on $\R^{3 \times 2N}$ can be written concisely as $g \phi(S)$ we define
\begin{equation}\label{eq:orthoact}
g(S) = \phi^{-1}\big(g\phi(S)\big)
\end{equation}
It is almost immediate from the definition that if $S \in D_N^6$ also $g(S) \in D_N^6$.

\begin{remark}\label{rem:notationb}
If $x$ and $y$ are different coordinates on phase space indicating the same point, we write $F(x)$ and $F(y)$ for the value of the function $F$ in that point although functional dependence may be different. Even if $S$ is a coordinate on a set of equivalence classes like the 'Albouy coordinates' we still write $F(S)$.
\end{remark}

%--------------------------------------------------------------------------------
\subsection{Clusters}\label{sec:clusters}
In section \ref{sec:albco} we considered sequences of states converging to a critical point or a critical point at infinity. In case of the latter we have in fact a critical orbit in mind in a three body problem, consisting for example of a single body and at a large distance two bodies rotating around a common centre of mass. Thus we have two \emph{clusters}, one of a single body and another of two bodies. We formalize this idea in the following definition of the $D_N$-cluster decomposition.
\begin{definition}[Cluster decomposition]\label{def:dndecomp}
Let $X = (\xi_1,\ldots,\xi_N) \in D_N$ and define the \emph{$D_N$-cluster decomposition}, into two clusters of $l$ and $N-l$ bodies, as $X = X_{\alpha} + X_{\beta} + X_{\alpha,\beta}$ with $X_{\alpha}, X_{\beta}, X_{\alpha,\beta} \in D_N$, where
\begin{align*}
X_{\alpha} &= (\alpha_1,\ldots,\alpha_l,0,\ldots,0)\\
X_{\beta} &= (0,\ldots,0,\beta_1,\ldots,\beta_{N-l})\\
X_{\alpha,\beta} &= (\alpha_0,\ldots,\alpha_0,\beta_0,\ldots,\beta_0)
\end{align*}
and the entries of $X_{\alpha}$ and $X_{\beta}$ are $\alpha_i = \xi_i - \alpha_0$ for $i=1,\ldots,l$ and $\beta_j = \xi_{l+j} - \beta_0$ for $j = 1,\ldots,N-l$. Furthermore $\alpha_0 = M_{\alpha}^{-1}\sum_{i=1}^l m_i \xi_i$ with $M_{\alpha} = \sum_{i=1}^l m_i$ and $\beta_0 = M_{\beta}^{-1}\sum_{i=l+1}^N m_i \xi_i$ with $M_{\beta} = \sum_{i=l+1}^N m_i$.

We say that $X$ is split into the two \emph{clusters} $X_{\alpha}$, $X_{\beta}$ and the \emph{centres of mass} $X_{\alpha,\beta}$. The projections $X \mapsto X_{\alpha}$ and $X \mapsto X_{\beta}$ are called the \emph{projections onto the clusters} and $X \mapsto X_{\alpha,\beta}$ is called the \emph{projection onto the centres of mass}.
\end{definition}
We use this definition of $D_N$-cluster decomposition to define a decomposition into clusters of a state $S$.
\begin{definition}
The \emph{decomposition into clusters} of the state $S = (X, Y, Z, P, Q, R)$ is obtained by applying the $D_N$-cluster decomposition to each component of $S$.
\end{definition}
Since the decomposition is in $D_N$, the construction can be iterated to a decomposition into $m$ clusters as long as $m \leq N$. A nice property in the next lemma is that the decomposition is orthogonal, which follows from a direct computation using definition \ref{def:dndecomp}.
\begin{lemma}\label{lem:orthodecomp}
The decomposition $X = X_{\alpha} + X_{\beta} + X_{\alpha,\beta}$ is orthogonal with respect to the inner product on $D_N$, that is $\inprod{X_{\alpha}}{X_{\beta}} = 0$, $\inprod{X_{\beta}}{X_{\alpha,\beta}} = 0$ and $\inprod{X_{\alpha,\beta}}{X_{\alpha}} = 0$.
\end{lemma}

Lemma \ref{lem:orthodecomp} together with lemma \ref{lem:functions} shows that the functions K, I and L are additive over clusters. This is useful in order to view each cluster as an independent $l$-body system. The latter is in however only useful when the distances between the centers of mass of the different clusters go to infinity in which case the potential $V$ also becomes asymptotically additive. For a precise statement, see lemma \ref{lem:additivity} in section \ref{sec:techs}. In view of these arguments we define the cluster separation property. Let $\Pi$ be the projection from phase space to configuration space. The following definition relates to Albouy coordinates.
\begin{definition}[Cluster separation property]\label{def:csp}
Let $N \geq 1$ and let $l$ be an integer with $1 \leq l \leq N$. A sequence $S_k \in D_N^6$ has the \emph{cluster separation property} if $l$ clusters of $S_k$ exists such that for $k \to \infty$:
\begin{enumerate}[i), itemsep=0.5ex, topsep=0pt]
\item each cluster in $\Pi(S_k)$ converges to a point,
\item for each pair in the cluster decomposition of $\Pi(S_k)$, the distance between the centres of mass converges to infinity.
\end{enumerate}
\end{definition}
It is not hard to prove that such sequences actually exist. Or more precisely that every sequence $S_k$ contains a subsequence having this property. Using the notion of (asymptotic) additivity of the functions $K$, $I$, $L$ and $V$, it makes sense to talk about clusters as smaller $l$-body problems in an $N$-body problem, of course when $l < N$. We will come back to this in section \ref{sec:techs}.
\begin{lemma}\label{lem:csp}
Every sequence $S_k$ has a subsequence with the cluster separation property.
\end{lemma}
\begin{proof}
We provide a construction in several steps. 1) A cluster containing a single body has the property that, restricted to this cluster, $\Pi(S_k)$ converges to a point. 2) Now suppose, possibly after renumbering the bodies, the cluster of the first $m$ bodies again has the property that, restricted to this cluster, $\Pi(S_k)$ converges to a point. For each of the remaining $N-m$ bodies, take the liminf of the distance this body to the centre of mass of the first cluster. If the liminf is finite add the body to the first cluster. If the liminf is infinite, add the body to the 'remaining bodies'. 3) We are left with at most $N-m$ remaining bodies. Now start again at step (1) with with these $N-m$ bodies.

The following remarks conclude the proof. One has to take a refined subsequence of $S_k$ at each step a body is added to a cluster. Since at each cycle through the steps above $m$ is at least one, the procedure is finite. By construction, condition ii) in definition \ref{def:csp} is satisfied. There are two extremes: all bodies in one cluster and every body in a separate cluster.
\end{proof}

In the sequel we will often encounter sequences decomposed into clusters. Unless stated otherwise we will always assume that such a decomposition of a sequence $S_k$ is independent of the sequence index $k$.

We end this section with rather obvious properties of sequences and clusters, but it is useful enough to warrant explicit mention. Roughly speaking if a critical sequence has a cluster decomposition which is constant over the sequence elements then each sequence of clusters is critical and vice versa. The next lemma provides a more precise formulation. The proof of the lemma follows almost immediately from the definition (\ref{def:critseq}) of a critical sequence.
\begin{lemma}\label{lem:critclusters}
Let $S_k$ be a sequence of states with a decomposition into $l$ clusters.
\begin{enumerate}[i), itemsep=0.5ex, topsep=0pt]
\item Suppose $S_k$ is a critical sequence with a compatible sequence of multipliers $\lambda_k$. Then the projection of $S_k$ on each of the $l$ sequences of clusters is a critical sequence.
\item Suppose that each of the $l$ sequences of clusters is critical. Then $S_k$ is a critical sequence if a common compatible sequence of multipliers exists, for which each sequence of clusters is critical.
\end{enumerate}
\end{lemma}

Furthermore the cluster decomposition commutes with the action of the rotation group defined in equation \eqref{eq:orthoact}. Again, this is almost immediately clear from the definition.

%--------------------------------------------------------------------------------
\section{Main results}\label{sec:cpimain}
The main theorem lists all possibilities for critical sequences.

\begin{theorem}[Main theorem on critical sequences]\label{the:mcs}
Let $S_k$ be a critical sequence of $H$.
\begin{enumerate}[i), itemsep=0.5ex, topsep=0pt]
\item \emph{Critical point at infinity.} Let $S_k$ be a horizontal critical sequence. Suppose at least two clusters exist in the sense of the cluster separation property, then the projection of the subsequence on each cluster converges to a (point on a) relative equilibrium.
\item \emph{Critical point or collision.} If the projection of $S_k$ on configuration space converges to a point, then a subsequence exists converging to a collision or a (point on a) relative equilibrium.
\end{enumerate}
\end{theorem}

The theorem says that by constructing critical sequences we not only get critical points at infinity but also 'ordinary' critical points and collisions. To prove the theorem above we again follow the approach of Albouy \cite{alb1993} adapted to our case.

The main theorem summarises a number of results that we now prove in a series of propositions each addressing one of the possibilities. Proposition \ref{pro:critseqrel} deals with the case of critical points outside the collision set. The critical sequences satisfying the conditions of proposition \ref{pro:critseq} converge to critical points of the Hamiltonian restricted to the level set of angular momentum, in fact the 'ordinary' critical points for which we do not need sequences at all. Proposition \ref{pro:horcritseq} considers horizontal critical sequences.

\begin{proposition}\label{pro:critseqrel}
Let $S_k$ be a critical sequence such that the projection on configuration space converges to a point outside the collision set. Then a subsequence exists that converges to a (point on a) relative equilibrium.
\end{proposition}

If we drop the condition on convergence to a point outside the collision set then the limit is either a relative equilibrium or a collision. Thus we see that in defining critical sequences we did not introduce that many new phenomena.

\begin{proposition}\label{pro:critseq}
Let $S_k$ be a critical sequence such that the projection on configuration space converges to a point. Then a subsequence exists converging either to a collision or to a (point on a) relative equilibrium.
\end{proposition}

Finally only the horizontal critical sequences yield critical points at infinity. The condition that angular momentum is constant on a horizontal critical sequence ensures that the only possible limit points of the sequence are relative equilibria, possibly in more than one cluster. If there are two or more clusters they move infinitely far apart along the sequence.

\begin{proposition}\label{pro:horcritseq}
Let $S_k$ be a horizontal critical sequence. Suppose at least two clusters exist in the sense of the cluster separation property, then the projection of the subsequence on each cluster converges to a (point on a) relative equilibrium.
\end{proposition}

%--------------------------------------------------------------------------------
\section{Proving the main results}\label{sec:promain}
In the proofs of propositions \ref{pro:critseqrel}, \ref{pro:critseq} and \ref{pro:horcritseq} we use a number of rather technical results which we provide in section \ref{sec:techs}.

There are a few notions and notations that we need in the proofs of the propositions. We review them here, more details are given in section \ref{sec:techs}. For a certain choice of coordinates, see lemma \ref{lem:mulco}, the expression $\Lgr(S, \lambda) = \grad H - \sum_i \lambda_i\, \grad L_i$ reduces to $\grad\,H - \norm{\lambda} \grad\,L_z$, these are called \emph{multiplier coordinates}. Using these coordinates we define $K_z = \norm{P}^2 + \norm{Q}^2$. A \emph{small sequence} is a critical sequence such that $I_k \to 0$, see definition \ref{def:smallseq}.

%--------------------------------------------------------------------------------
\subsection{Proofs of propositions}
Boundedness is the key ingredient in the proof of proposition \ref{pro:critseqrel}, but also the fact that the moment of inertia is bounded away from zero. Roughly speaking, this finiteness allows us to conclude that the limit is a (point on a) relative equilibrium.

\begin{proof}[Proof of proposition \ref{pro:critseqrel}]
$I_k$ is bounded, because the projection of the sequence $S_k$ on configuration space converges to a point. Moreover, since this point is not a collision, $V_k$ is bounded and $I_k$ is bounded away from zero. From lemma \ref{lem:kplusv} we see that also $K_k$ is bounded, thus we conclude that $\norm{S_k}^2 = \norm{I_k} + \norm{K_k}$ is bounded. Therefore the sequence $S_k$ is bounded. Thus $S_k$ has a limit point and hence a convergent subsequence. From now on we restrict to the latter. Again using lemma \ref{lem:kplusv} we know that $\norm{K_k}$ and $\norm{V_k}$ have the same limit and if we use \emph{multiplier coordinates}, $K_{z,k}$ and $\norm{V_k}$ have the same limit. Now suppose that $V$ is non-zero, then $K_{z,k}$ is bounded and bounded away from zero. From lemma \ref{lem:estims} iii) we see that $L_{z,k}$ is also bounded away from zero. Thus the limit satisfies the equation for a critical point of $H$ restricted to a non-zero level of $L$, and the subsequence converges to a (point on a) relative equilibrium. It remains to prove that $V_k$ is non-zero, or directly that $K_k$ is bounded away from zero. In the case of gravitation $V$ can only be zero if the bodies go infinitly far apart, which they do not by the assumptions of the proposition. In the Coulomb case however, $V$ can be zero for bodies at finite distances. We now use the fact that we do not consider the level zero of angular momentum. Then $K$ must be non-zero.
\end{proof}

The key ingredient in the next proof is the compatibility of the multipliers: a sequence converging to a collision has multipliers not compatible with a sequence converging to a relative equilibrium and vice versa.

\begin{proof}[Proof of proposition \ref{pro:critseq}]
Suppose the critical sequence is such that the projection on configuration space converges to a collision point (non-collision is the case of proposition \ref{pro:critseqrel}). Let us decompose the sequence into clusters of bodies having the same limit when projected on configuration space. Then on each cluster the projection of $I$ converges to zero and using lemma \ref{lem:smallseq} the sequence of compatible multipliers converges to infinity. However the projection on the centres of mass is not a collision and thus a subsequence converging to a relative equilibrium exists and the accompanying subsequence of compatible multipliers has a finite limit. This contradicts lemma \ref{lem:critclusters} and thus a critical sequence either: 1) converges to a collision; or 2) a subsequence converging to a (point on a) relative equilibrium exists.
\end{proof}

The proof of the last proposition hinges on the fact that for a sequence with the cluster decomposition property, the sequences of projections on the clusters are independent with respect to kinetic and potential energy.

\begin{proof}[Proof of proposition \ref{pro:horcritseq}]
According to lemma \ref{lem:csp}, we may decompose the sequence (or a subsequence) into clusters. For each cluster we get a sequence by projecting the original sequence on the current component of the cluster decomposition. Now according to proposition \ref{pro:critseq}, each of these sequences either is a small sequence or has a subsequence converging to a (point on a) relative equilibrium. The possibility of a small sequence can be ruled out because $H \to -\infty$. This follows from the fact that a small sequence is a critical sequence and from lemma \ref{lem:kplusv} we have $K + V \to 0$, but on a small sequence $V \to -\infty$ and therefore $H = \frac{1}{2} K + V = \frac{1}{2} (K + V) + \frac{1}{2} V \to \frac{1}{2} V \to -\infty$, also see corollary \ref{cor:nocol}. Thus we arrive at a contradiction. We are left with sequences like in proposition \ref{pro:critseqrel}, converging to a (point on a) relative equilibrium in a subsystem. The sequence of multipliers has a finite (non-zero) limit. The centres of mass of the clusters form a critical sequence with the same sequence of multipliers. Since the centres go infinitely far apart, the potential approaches zero and therefore this sequence is a critical sequence of $K$. Now from lemmas \ref{lem:kcrita} and \ref{lem:kcritb} we infer that $K \to 0$ because on the horizontal sequence we started with, angular momentum has a fixed, finite value.
\end{proof}

%--------------------------------------------------------------------------------
\subsection{Some technical results}\label{sec:techs}
Here we collect a number of technical results used in the proofs of the previous section.

Given a state or a sequence of states, we will at various instances use a partition into clusters. Many functions we need are additive over these clusters, namely the moment of inertia $I$, twice the kinetic energy $K$ and angular momentum $L$. One exception being the potential. However on a sequence with the cluster separation property, taking the limit the potential is additive over the clusters. This is made more precise in the following lemma. 

\begin{lemma}\label{lem:additivity}
Let $S_k$ be a sequence of states with the cluster separation property. Suppose $S_k$ is decomposed into $l$ clusters denoted by $S_{k,j}$ with $j \in \{1,\ldots,l\}$. Let $F$ be one of the functions $K$, $I$, $L$, $V$ or the (components of the) vector valued function $\grad\, V$. Then the following holds
\begin{equation}
\lim_{k \to \infty} F(S_k) = \sum_{j=1}^l \lim_{k \to \infty} F(S_{k,j}).
\end{equation}
\end{lemma}

\begin{proof}
For the functions $K$, $I$ and $L$ we may drop the limits and still have equality. Namely, let $\inprod{\cdot}{\cdot}$ be an inner product and $\norm{\cdot}$ be the associated norm. If $x$ and $y$ are orthogonal with respect to $\inprod{\cdot}{\cdot}$, then $\norm{x + y}^2 = \norm{x}^2 + \norm{y}^2$. The $D_N$-decomposition is orthogonal, see lemma \ref{lem:orthodecomp}, with respect to the inner product from definition \ref{def:dispo}. The potential function $V$ is not defined using the inner product and equality without the limits does not hold indeed. However, the contribution to the potential function arising from interaction between bodies from different clusters converges to zero for $k \to \infty$. The reason is that the distance between different clusters converges to infinity and the potential approaches zero when distances go to infinity. Therefore in the limit $k \to \infty$ the only non-zero contribution to the potential comes from within the clusters. The same holds for the gradient of the potential, see section \ref{sec:set}.
\end{proof}

Note that we have not excluded the possibility that the potential converges to infinity along a (sub)sequence, that is the lemma holds even in the case that one of the clusters converges to a collision. The next lemma relates kinetic and potential energy in the limit along a critical sequence.
\begin{lemma}\label{lem:kplusv}
Let $S_k$ be a critical sequence of $H$ and let $V$ be a function homogeneous of degree $-1$, then
\begin{equation*}
\lim_{k \to \infty} \frac{(K_k + V_k)^2}{I_k + K_k} = 0.
\end{equation*}
\end{lemma}

\begin{proof}
From the definition of a critical sequence, see definitions \ref{def:compmul} and \ref{def:critseq}, we know that a compatible sequence of multipliers exists such that $e_k = \Lgr(S_k, \lambda_k)$ converges to zero as $k \to \infty$. Let $f$ be a vector field of unit length, then also $\lim_{k \to \infty} \inprod{e_k}{f} = 0$. This holds in particular when $f$ is tangent to the levels of $L$. Let us take such a vector field, namely $f=(X, Y, Z, -P, -Q, -R)/\sqrt{I+K}$. Then we have (suppressing the subscript $k$)
\begin{align*}
\inprod{e}{f}^2 &= \inprod{\grad\,H}{f}^2\\
&= \Bigg(X \partieel{X}{V} + Y \partieel{Y}{V} + Z \partieel{Z}{V} - \tfrac{1}{2}\bigg(P \partieel{P}{K} + Q \partieel{Q}{K} + R \partieel{R}{K}\bigg)\Bigg)^2 \bigg/\bigg(I + K\bigg)\\
&= (V + K)^2 / (I + K).
\end{align*}
In the above we used $H = \tfrac{1}{2}K + V$ and the fact that $V$ is homogeneous of degree $-1$. Since $\lim_{k \to \infty} \inprod{e_k}{f} = 0$ we also have $\lim_{k \to \infty} (K + V)^2 / (I + K) = 0$. It is easily checked that 1) $f$ is a unit vector field since $\norm{(X, Y, Z, -P, -Q, -R)} = \sqrt{I+K}$ and 2) $f$ is tangent to the levels of $L$.
\end{proof}

\begin{remark}
Let us make a few remarks on this lemma.
\begin{enumerate}[i), itemsep=0.5ex, topsep=0pt]
\item The lemma is formulated for a homogeneous potential $V$ of degree $-1$. It can easily be generalized for a potential homogeneous of degree $-\mu$, with $\mu > 0$, then the result is $\lim_{k \to \infty} \frac{(K_k + \mu V_k)^2}{I_k + K_k} = 0$.
\item At first sight it might seem strange that the value of the potential matters: adding an arbitrary constant to the potential $V$ does not affect the equations of motion. However, we made essential use of the fact that $V$ is a homogeneous function, a property that is lost if we add a constant to $V$. Another reason not to add a constant to $V$ is that we want the interaction to vanish when bodies go infinitely far apart. Thus it seems natural to require that in that case the potential energy approaches zero as well.
\item We conjecture that homogeneity of the potential is not essential. But the potential and its gradient vanishing when bodies go infinitely apart does seem essential. It is less clear how important the behaviour of the potential near collisions really is, apart from a convenient way to distinguish a collision by the values of the potential. The proof of lemma \ref{lem:kplusv} rests on a particular choice of the vector field $f$, other choices may also give useful estimates. An example of a potential for which we would like to obtain similar results is the Lennard-Jones potential, see \cite{reichl2016} which is repelling for bodies at a short distance, but attracting for bodies far apart. This potential can be regarded as being 'asymptotically homogeneous'.
\end{enumerate}
\end{remark}
Since kinetic energy $\frac{1}{2} K$ is always non-negative, lemma \ref{lem:kplusv} implies the following property of critical sequences.
\begin{corollary}\label{cor:nocol}
For a homogeneous potential of degree $-\mu < 0$, a critical sequence, or a subsequence thereof, never approaches a collision of repelling bodies.
\end{corollary}

In the sequel we will frequently use \emph{multiplier coordinates} such that the equation for a critical sequence takes a simpler form, namely $\Lgr(S_k, \lambda_k) = \grad\,H_k - \norm{\lambda_k} \grad\,L_{z,k}$. For future reference we write the two terms in components:
\begin{equation}\label{eq:lagrangecmps}
\begin{aligned}
\grad\,H &= \big(\partieel{X}{V},\partieel{Y}{V}, \partieel{Z}{V}, P, Q, R \big)\\
\grad\,L_z &= \big(Q, -P, 0, -Y, X, 0 \big)
\end{aligned}
\end{equation}
The above is based on the next lemma.
\begin{lemma}\label{lem:mulco}
Let $S_k$ be a critical sequence with compatible multipliers $\lambda_k$ like in definition \ref{def:compmul}. For each $k$, there is a $g \in SO(3)$ such that in rotated coordinates, called \emph{multiplier coordinates}, $\Lgr(S_k, \lambda_k) = \grad\,H_k - \norm{\lambda_k} \grad\,L_{z,k}$. 
\end{lemma}
\begin{proof}
To simplify notation we suppress the index $k$. Upon changing coordinates $S \mapsto g(S)$, see equation \eqref{eq:orthoact}, the function $H$ is invariant, but $L \mapsto gL$. Now we write $\Lgr(S, \lambda) = \grad\,H - \sum_i\lambda_i \grad\,L_i$ as $\grad\,H - \grad\,\inprod{\lambda}{L}$. Now $\inprod{\lambda}{gL} = \inprod{g^t\lambda}{L}$ and by a suitable choice of $g$ we get $g^t\lambda = (0,0,\norm{\lambda})$. We should note that the gradient is not invariant upon changing coordinates, but since we are only interested in the value zero we can safely ignore this.
\end{proof}

Another application of lemma \ref{lem:mulco} is that we are able to obtain simple expressions for estimates on $R$, $K$, $I$ and $L$ which are independent of the potential in the system.

\begin{lemma}\label{lem:estims}
Consider a critical sequence $S_k$ together with its compatible sequence of multipliers $\lambda_k$. In a coordinate system of lemma \ref{lem:mulco} we define $I_z = \norm{X}^2 + \norm{Y}^2$ and $K_z = \norm{P}^2 + \norm{Q}^2$. Then we have the following estimates for $k \to \infty$ (for sake of readability we suppressed the index $k$):
\begin{enumerate}[i), itemsep=0.5ex, topsep=0pt]
\item $\norm{R} = \smallo(1)$
\item $\sqrt{K_z} = \norm{\lambda} \sqrt{I_z} + \smallo(1)$
\item $\sqrt{K_z} = L_z/\sqrt{I_z} + \smallo(1)$
\end{enumerate}
\end{lemma}

\begin{proof}
Using \emph{multiplier coordinates} we have $e_k = \Lgr(S_k, \lambda_k) = \grad\, H_k - \norm{\lambda_k} \grad\, L_{z,k}$ converging to zero as $k \to \infty$, since $S_k$ is a critical sequence. We write $e_k$ in components. However we will only use the last three and since the first three would need some explanation we skip them:
\begin{equation*}
e_k = (*, *, *, P + \norm{\lambda} Y, \;Q - \norm{\lambda} X, \;R).
\end{equation*}
\textit{ad i)} Project $e_k$ on the sixth direction, '$R$-space'. Because $\norm{e_k} \to 0$ when $k \to \infty$ we have $\norm{R} = \smallo(1)$.\\
\textit{ad ii)} Consider the projection of $e_k$ on the fourth and fifth directions, '$(P,Q)$-space'. Then we have $(P + \norm{\lambda} Y, Q - \norm{\lambda} X) = \smallo(1)$. Using the triangle inequality we have
\begin{equation*}
\smallo(1) = \norm{(P, Q) - \norm{\lambda} (-Y,X)} \geq \mid \sqrt{K_z} - \norm{\lambda} \sqrt{I_z} \mid
\end{equation*}
or $\sqrt{K_z} = \norm{\lambda} \sqrt{I_z} + \smallo(1)$.\\
\textit{ad iii)} Starting again with $(P + \norm{\lambda} Y, Q - \norm{\lambda} X) = \smallo(1)$, note that also the inner product with the unit vector $(Y,-X)/\sqrt{I_z}$ is $\smallo(1)$. Computing this inner product we get
\begin{align*}
\Big\langle(P, Q) - \norm{\lambda} (Y, -X), (-Y, X)\Big\rangle/\sqrt{I_z} &= \Big(\inprod{-P}{Y} + \inprod{Q}{X}\Big)/\sqrt{I_z} + \norm{\lambda} \sqrt{I_z}\\ &= L_z/\sqrt{I_z} + \norm{\lambda} \sqrt{I_z}.
\end{align*}
Therefore $L_z/\sqrt{I_z} + \norm{\lambda} \sqrt{I_z} = \smallo(1)$ and using ii) we arrive at $\sqrt{K_z} = L_z/\sqrt{I_z} + \smallo(1)$.
\end{proof}

The following lemma will be useful to distinguish critical sequences converging to a collision. But first we define the notion of a \emph{small sequence} which we will need later on as well.

\begin{definition}\label{def:smallseq}
A \emph{small sequence} is a critical sequence such that $I_k \to 0$ when $k \to \infty$.
\end{definition}

\begin{lemma}\label{lem:smallseq}
Suppose $S_k$ is a small sequence, then both $H_k$ and the norm of the multiplier become unbounded.
\end{lemma}

\begin{proof}
From lemma \ref{lem:kplusv} we know that $K_k + V_k \to 0$ for $k \to \infty$. Furthermore $V_k$ becomes unbounded when $I_k$ converges to zero and therefore also $K_k$ becomes unbounded. Since $H_k = \frac{1}{2} K_k + V_k = \frac{1}{2}(K_k + V_k) + \frac{1}{2} V_k \to \frac{1}{2} V_k$, $H_k$ becomes unbounded. Using \emph{multiplier coordinates} we see from lemma \ref{lem:estims} ii) that $\norm{\lambda}$ too becomes unbounded.
\end{proof}

In view of corollary \ref{cor:nocol} we have the following corollary of lemma \ref{lem:smallseq}
\begin{corollary}\label{cor:nocolh}
Since a critical sequence or a subsequence thereof will never approach a collision of repelling bodies, also in the case of the Coulomb potential we have $H_k \to -\infty$ along a small sequence.
\end{corollary}

\commentaar{
\begin{remark}
When we apply the methods in the proof of lemma \ref{lem:kplusv} to a potential $V$ homogeneous of degree $-\mu$ we get $K_k + \mu V_k \to 0$ for $k \to \infty$. Then $H_k = \frac{1}{2} K_k + V_k = \frac{1}{2}(K_k + \mu V_k) + (1-\frac{\mu}{2}) V_k \to (1-\frac{\mu}{2}) V_k$, which indicates that $\mu = 2$ might be an exceptional case.
\end{remark}
}

\commentaar{toelichting of waarom}

On a horizontal critical sequence (see definition \ref{def:horcritseq}) the kinetic energy converges to zero. This follows from the following lemma on critical sequences of $K$.

\commentaar{free bodies, here centres of mass of clusters}

\begin{lemma}\label{lem:kcrita}
Consider a critical sequence of $K$ such that $K$ is bounded away from zero. Then $\norm{L}$ converges to infinity.
\end{lemma}

\begin{proof}
Again we use \emph{multiplier coordinates}. From the multiplier equations, see equation \eqref{eq:lagrangecmps}, we have: $\norm{\lambda} \sqrt{K_z} = \smallo(1)$. On the other hand from the estimates that do not depend on the potential, lemma \ref{lem:estims} ii), we have $\sqrt{K_z} = \norm{\lambda} \sqrt{I_z} + \smallo(1)$. Then we have two expressions for $\sqrt{K_z}$, yielding $K_z = \smallo(\sqrt{K_z}) + \smallo(\sqrt{I_z})$, implying $I_z$ becoming unbounded. Combining the latter with estimate \ref{lem:estims} iii) namely $\sqrt{K_z} = L_z/I_z + \smallo(1)$, yields that $L_z$ becomes unbounded.
\end{proof}

\begin{lemma}\label{lem:kcritb}
Consider a critical sequence of $K$ such that $L_z$ is bounded away from zero. Then $K$ and $I_z$ converge to zero.
\end{lemma}

\begin{proof}
Like in the proof of the previous lemma, we use \emph{multiplier coordinates}. From the estimates lemma \ref{lem:estims} i) $R = \smallo(1)$ and $\norm{\lambda} \sqrt{K_z} = \smallo(1)$ in the proof of the previous lemma, we see that $K$ converges to zero. Combining the estimates in lemma \ref{lem:estims} ii) $\sqrt{K_z} = \norm{\lambda} \sqrt{I_z} + \smallo(1)$ and $\norm{\lambda} \sqrt{K_z} = \smallo(1)$ we get two expressions for $\norm{\lambda} \sqrt{K_z}$, yielding $\norm{\lambda}^2 \sqrt{I_z} = \smallo(1) + \smallo(\norm{\lambda})$. From the latter we infer that $I_z$ converges to zero.
\end{proof}

%--------------------------------------------------------------------------------
\section{Critical points at infinity in the charged three body problem}\label{sec:cpictbp}
Now that we know that the critical points at infinity can be found from the so called horizontal critical sequences, we wish to construct these for the charged three body problem. In particular we need to construct a horizontal critical sequence with the cluster separation property. There are only three possibilities to form clusters, namely 1) three clusters each containing only one body; 2) two clusters of one and two bodies respectively, and 3) one cluster containing all three bodies. We can immediately rule out the last case because a sequence with only a single cluster in the sense of the cluster separation property, can not converge to a critical point at infinity. We will not consider the first case because then angular momentum $L$ would have to converge to zero, but we start with a horizontal critical sequence on a non-zero level of $L$. Thus we are left with the case of a horizontal critical sequence with two clusters. As it turns out there is only a limited number of possibilities, for each we compute the value of the bifurcation parameter $\nu = -h \norm{l}^2$, where $h$ is the value of $H$, the energy, and $l$ is the value of angular momentum $L$.

Below we consider a horizontal critical sequence with two clusters using multiplier coordinates. The fact that $I_z$ tends to zero on a horizontal sequence (see lemma \ref{lem:kcritb}), implies that the two clusters go infinitely apart in the $z$ direction. This in turn implies that there are no critical points at infinity in planar three body problems with gravitational or Coulomb potentials.

\textbf{Construction of a horizontal critical sequence} We will in particular construct a horizontal critical sequence $S_k$ consisting of two clusters having the cluster decomposition property. Like in definition \ref{def:dndecomp} we call the components of the clusters $X_{\alpha}$ and $X_{\beta}$ and similarly for the other components. In the cluster called $\alpha$ we put two bodies with masses $m_1$ and $m_2$ and in the cluster called $\beta$ we put one body with mass $m_3$. Since we have only two clusters, the centres of mass lie on a line which we take to be the $z$-axis. Thus we have $X_{\alpha,\beta} = (0,0,0)$, $Y_{\alpha,\beta} = (0,0,0)$ and $Z_{\alpha,\beta} = (\frac{m_3}{m_1+m_2} z_k, \frac{m_3}{m_1+m_2} z_k, -z_k)$, where $\lim_{k \to \infty} z_k = \infty$. The single body in the cluster $\beta$ just remains at its centre of mass. So we have $X_{\beta} = Y_{\beta} = Z_{\beta} = (0,0,0)$. The two bodies in cluster $\alpha$ form a two body problem which for a $SO(3)$-invariant potential is planar and in view of the remarks preceding this paragraph, we take this plane to be the '$(x,y)$-plane'. Therefore we take $X_{\alpha} = (\xi_1, \xi_2, 0)$, $Y_{\alpha} = (\eta_1, \eta_2, 0)$ and $Z_{\alpha} = (0,0,0)$ such that $m_1 \xi_1 + m_2 \xi_2 = 0$ and similarly for $\eta$. So far we only have clusters. Projecting the sequence $S_k$ on the $\beta$ component we get a body at rest at infinity when taking the limit $k \to \infty$. But projecting on the $\alpha$ component the limit should be a relative equilibrium of the two body problem. 

\textbf{The planar two body problem} We consider the planar two body problem with potential $V$ on the '$(x,y)$-plane' of the previous paragraph. That is, we consider a Hamiltonian system on $T^*(\R^4)$ with standard symplectic form and coordinates $x_i = (\xi_i, \eta_i)$ and corresponding momenta $y_i$ with Hamiltonian
\begin{equation*}
H(x, y) = \tfrac{1}{2} \bigg( \frac{\norm{y_1}^2}{m_1} + \frac{\norm{y_2}^2}{m_2} \bigg) + V(x_1, x_2)
\end{equation*}
where the potential $V$ depends on relative distance only $V(x_1, x_2) = \gamma f(\norm{x_1 - x_2})$ where $\gamma$ depends only on properties of the bodies (like mass and charge) and $f$ is a function. The function $f$ will eventually be defined as $f(x) = \frac{1}{x}$. But for the moment we leave it unspecified. It is easily seen that this Hamiltonian system is translation and rotation invariant. To fulfill the condition of the cluster, namely that the centre of mass remains at zero, we apply a (symplectic) change of coordinates to a centre of mass frame.
\begin{lemma}\label{lem:comf}
By the symplectic change of coordinates
\begin{equation*}
%(x_1, x_2, y_1, y_2) \mapsto (q_1, q_2, p_1, p_2) = (x_1 - x_2, \frac{m_1 x_1 + m_2 x_2}{m_1 + m_2}, \frac{m_2 y_1 - m_1 y_2}{m_1 + m_2}, y_1 + y_2)
(x_1, x_2, y_1, y_2) \mapsto (q_1, p_1, q_2, p_2) = (x_1 - x_2, \frac{m_2 y_1 - m_1 y_2}{m_1 + m_2}, \frac{m_1 x_1 + m_2 x_2}{m_1 + m_2}, y_1 + y_2)
\end{equation*}
the new Hamiltonian system is separated on $T^*(\R^2) \times T^*(\R^2)$ with standard symplectic form on each component we get
\begin{align*}
H(q, p) &= \tfrac{1}{2} \bigg( \Big(\frac{1}{m_1} + \frac{1}{m_2} \Big) \norm{p_1}^2 + \frac{1}{m_1 + m_2} \norm{p_2}^2 \bigg) + \gamma f(\norm{q_1})\\
        &= \frac{1}{2\mu} \norm{p_1}^2 + \gamma f(\norm{q_1}) + \frac{1}{2m} \norm{p_2}^2\\
        &= H_1(q_1, p_1) + H_2(q_2, p_2).
\end{align*}
Here $\mu = \frac{m_1 m_2}{m_1 + m_2}$ is the reduced mass and $m = m_1 + m_2$. The second Hamiltonian $H_2$ describes the motion of the centre of mass. By setting $q_2 = p_2 = 0$ the centre of mass remains at rest.
\end{lemma}
The proof of the lemma is a straightforward calculation. We now concentrate on the first Hamiltonian system on $\R^4$ with coordinates $(q_1, p_1)$, dropping the no longer necessary subscript.

The system with Hamiltonian $H(q, p) = \frac{1}{2\mu} \norm{p}^2 + \gamma f(\norm{q})$ is still rotation invariant: angular momentum is a conserved function. Recall that we are looking for a relative equilibrium, in other words a critical point of $H$ restricted to a level of angular momentum. One way to solve this is to use the Lagrange multiplier method and recognize a problem of finding central configurations for this problem, see \cite{hwz2019}. Another way is to reduce with respect to the symmetry group and find stationary points for the reduced system. We will adopt the last approach. To proceed we introduce polar coordinates in the $(q_1, q_2)$-plane: $q_1 = r \cos \phi$ and $q_2 = r \sin \phi$. (Note that $q_1, q_2 \in \R$ refer to the components of $q$, not to be confused with $q_1$ and $q_2$ in the lemma.) The result is the following, see \cite{arnold1980} or \cite{am1987}.

\begin{proposition}\label{pro:tbpreduc}
The reduction of the two body problem is a Hamiltonian system on $T^*(\R_{>0})$ with standard symplectic form and Hamiltonian
\begin{equation*}
H(r, p_r) = \frac{1}{2\mu} p_r^2 + U(r), 
\end{equation*}
where $U(r) = \frac{\ell^2}{2 \mu r^2} + \gamma f(r)$ and $\ell$ is the value of angular momentum.
\end{proposition}

The reduced system has a critical point at $(r, p_r) = (r^*, 0)$ where $r^*$ is a solution of $U'(r) = 0$. Taking $f(x) = \frac{1}{x}$ we find a single solution $r^* = -\frac{\ell^2}{\mu \gamma}$, provided that $\gamma < 0$. If $\gamma > 0$ then the reduced system does not have a critical point. Thus the two body problem with $f(x) = \frac{1}{x}$ has a relative equilibrium if $\gamma < 0$.

\textbf{Critical values related to critical points at infinity} Now that we have a horizontal critical sequence with the cluster separation property we are able to compute the values of the bifurcation parameter $\nu = -h \ell^2$ where $h$ is the value of $H$ and $\ell$ is the value of $L$. From proposition \ref{pro:horcritseq} and lemmas \ref{lem:kcrita}, \ref{lem:kcritb} we know that the contribution of the projection on the centres of mass to the energy $H$ is zero for $k \to \infty$. Furthermore the projection on the cluster containing one single body is a small sequence, \ref{lem:smallseq} so the contribution to kinetic energy is zero in the limit. Therefore the only contribution to $H$ is from the projection on the clusters containing two bodies. To proceed we have to make a distinction between the gravitation potential and the Coulomb potential.

In case of the gravitation potential the constant $\gamma$ equals $-m_1 m_2$ which is always negative. Therefore there is always a unique critical point of the reduced Hamiltonian and thus a unique relative equilibrium in the two body problem. The value of the energy is $h = H(r^*, 0) = H(-\frac{\ell^2}{\mu \gamma}, 0) = - \frac{\mu \gamma^2}{2 \ell^2}$ and thus $\nu = -h\ell^2 = \frac{1}{2}\mu \gamma^2 = \frac{(m_1 m_2)^3}{2(m_1 + m_2)}$. Since there are three combinations of two masses out of three, there are three critical values related to critical points at infinity:
\begin{equation*}
\nu_3 = \frac{(m_1 m_2)^3}{2(m_1 + m_2)}, \;\;\; \nu_1 = \frac{(m_2 m_3)^3}{2(m_2 + m_3)}, \;\;\; \nu_2 = \frac{(m_1 m_3)^3}{2(m_1 + m_3)}.
\end{equation*}

The Coulomb force is only attracting if the charges have different sign. The constant $\gamma$ in this case equals $c_i c_j$ where $c_i$ is the electrical charge of body $i$. Given three bodies with differently signed charges, two signs are always equal. Let us assume the charges of bodies 1 and 2 have equal signs. Then only two combinations form a cluster converging to a relative equilibrium at infinity. Thus there are two critical values related to critical points at infinity:
\begin{equation*}
\nu_1 = \frac{m_1 m_3 (c_1 c_3)^2}{2(m_1 + m_3)}, \;\;\; \nu_2 = \frac{m_2 m_3 (c_2 c_2)^2}{2(m_2 + m_3)}.
\end{equation*}

%--------------------------------------------------------------------------------
%\section{Discussion and comments}\label{sec:disco}
\commentaar{moet dat?}
\commentaar{kritieke punten op oneindig niet in familie van eerder gevonden kritieke punten, ihb niet onder dilatie symmetrie}
\commentaar{geen kritieke puten op oneindig in een vlak systeem, zie boven}

%----------------------------------------------------------------------
\bibliographystyle{plain}
\bibliography{cpibib}

\end{document}